\documentclass[12pt]{amsart}

\usepackage{graphicx}
\usepackage{amsfonts}
\usepackage{amsmath}
\usepackage{amssymb}
\usepackage{amscd, amsthm}
\usepackage{verbatim}
\usepackage{latexsym}
\newcommand{\ord}{\mathrm{ord}}

\newcommand{\cG}{\mathcal{G}}
\newcommand{\cO}{\mathcal{O}}

\newcommand{\A}{\mathbb{A}}
\newcommand{\Q}{\mathbb{Q}}
\newcommand{\R}{\mathbb{R}}

\newcommand{\C}{\mathbb{C}}

\newcommand{\Z}{\mathbb{Z}}

\newcommand{\aZl}{\overline{\mathbb{Z}_{\ell}}}
\newcommand{\aQl}{\overline{\mathbb{Q}_{\ell}}}

\newcommand{\Aut}{\mathrm{Aut}\,}
\newcommand{\GL}{\mathrm{GL}}
\newcommand{\PGL}{\mathrm{PGL}}
\newcommand{\PSL}{\mathrm{PSL}}
\newcommand{\GSp}{\mathrm{GSp}}
\newcommand{\PSp}{\mathrm{PSp}}
\newcommand{\PGSp}{\mathrm{PGSp}}

\newcommand{\Sp}{\mathrm{Sp}}
\newcommand{\SO}{\mathrm{SO}}

\newcommand{\Spin}{\mathrm{Spin}}
\newcommand{\SU}{\mathrm{SU}}

\newcommand{\rhobar}{\overline{\rho}}

\newtheorem{theorem}{Theorem}[section]
\newtheorem{lemma}[theorem]{Lemma}
\newtheorem{prop}[theorem]{Proposition}

\newtheorem{cor}[theorem]{Corollary}

\def\rhobar{ {\bar {\rho} } }

\newcommand{\F}{{\mathbb F}}
\newcommand{\SL}{\mathrm{SL}}

\newcommand{\Ind}{\mathrm{Ind}}

\begin{document}

\title[Functoriality and the  Inverse Galois problem]
{Functoriality and the Inverse Galois problem} 
%\\ \small{\tt preliminary version}}
\author[C.~Khare]{Chandrashekhar Khare}
\email{shekhar@math.utah.edu}
\address{Department of Mathematics \\
         University of Utah \\
         155 South 1400 East, Room 233 \\
           Salt Lake City, UT 84112-0090 \\
         U.S.A.}\thanks{CK was partially supported by NSF grant DMS 0355528.}
\author[M.~Larsen]{Michael Larsen}
\email{larsen@math.indiana.edu}
\address{Department of Mathematics\\
     Indiana University \\
     Bloomington, IN 47405\\
     U.S.A.}
 \author[G.~Savin]{Gordan Savin}
\email{savin@math.utah.edu}
\address{Department of Mathematics \\
         University of Utah \\
         155 South 1400 East, Room 233 \\
           Salt Lake City, UT 84112-0090 \\
         U.S.A.}\thanks{GS was partially supported by NSF grant DMS 0551846.}

\date{}

\begin{abstract}
We prove that for any prime $\ell$ and any even integer $n$, there are infinitely many exponents $k$ for which $\PSp_n(\F_{\ell^k})$ or $\PGSp_n(\F_{\ell^k})$  appears as a Galois group over $\Q$. This generalizes a result of Wiese  \cite{Wiese}, which inspired this paper. 
\end{abstract}

\maketitle

\section{Introduction}

The inverse Galois problem asserts that every finite group $G$ occurs as 
${\rm Gal}(K/\Q)$ for $K/\Q$ a finite Galois extension of $\Q$. This has received much attention. It is natural to focus first on simple groups $G$.
The first infinite family of non-abelian finite simple groups for which the problem was solved was the family of alternating groups.  Hilbert proved his irreducibility theorem for this purpose, thus showing that it suffices to prove that $A_n$ occurs as the Galois group of a finite regular extension of $\Q(T)$.

The main advance on this problem in recent decades is the rigidity 
technique of Thompson \cite{Thompson}. This method has solved the problem for most of the  sporadic groups: it realizes all sporadic groups with the exception of  the Mathieu groups  $M_{23}$ and $M_{24}$ as Galois groups of regular extensions of $\Q(T)$.
We refer to \cite{det}, and the references therein, for results  towards the inverse Galois problem that are proved by the rigidity method and its variants.

For classical groups, rigidity-type methods have met with only sporadic success. Typically  these methods   seem
to work for $G(\F_{\ell^k})$, with $G$ a Chevalley group over the prime field  $\F_\ell$, when  $k$ is \textit{small} as compared to the rank of $G$.

Recently, Wiese \cite{Wiese} proved a result of the opposite kind:

\begin{theorem}
\label{inspiration}
Let $\ell$ be any prime.  Then there exist infinitely many integers $k$ such that
at least one of $\PSL_2(\F_{\ell^k})$ and $\PGL_2(\F_{\ell^k})$ can be realized as
a Galois group over $\Q$.  In particular, there are infinitely many integers $k$ for which the finite
simple group $L_2(2^k) = \PSL_2(\F_{2^k}) = \PGL_2(\F_{2^k})$ can be realized.
\end{theorem}

This paper generalizes Wiese's result to finite simple groups of symplectic type.

\begin{theorem}\label{main}
 If we fix a prime $\ell$ and integers $n,t \geq 1$ 
 with $n=2m$ even, the finite group $\PSp_{n}(\F_{\ell^k})$ or $\PGSp_n(\F_{\ell^k})$ occurs as
 a Galois group over $\Q$ for some integer $k$ divisible by $t$.
\end{theorem}

The method of \cite{Wiese} relies on results in \cite{KW}. In particular it relies on \cite[Lemma 6.3]{KW},  which asserts that if one ensures certain ramification properties of a compatible system of $2$-dimensional representations of $G_\Q$, then its residual representations for small residue characteristics are large.  Wiese uses this lemma and some other techniques  and results from \cite{KW}.
One may remark, however, that given some constructions of automorphic
forms, the only result from \cite{KW} one really needs 
to use is the simple but crucial \cite[Lemma 6.3]{KW} .

To prove our  theorem we construct a continuous irreducible 
representation $\rho: G_\Q \rightarrow \GL_n(\aQl)$
that is unramified outside $\ell$, the infinite place $\infty$, and another auxiliary prime $q$, and whose image is contained in either the orthogonal or symplectic similitudes.   
The representation $\rho$ is constructed so that the image of $\rho(D_q)$, with $D_q$ a decomposition group  at $q$ in $G_\Q$, is a metacyclic group,
which acts irreducibly on $\aQl^n$ and preserves an alternating form, up to
a multiplier.  Thus one knows that the image of $\rho$ is contained in fact in the symplectic similitudes. We ensure  that the order of $\rho(I_q)$, with $I_q$ an inertia group at $q$ of $G_\Q$,  is a prime $p \neq \ell$ that is larger than some large integer $N$. The representation $\rho$ has the property that all open  subgroups  $H$ of index $\leq N$  contain the image of $\rho(D_q)$. (The $N$ here is the $N_d$ of Theorem \ref{crucial}.)   This is ensured by choosing $q$ to split in all extensions of $\Q$ of degree at most $N$ that  are unramified outside $\ell$ and $\infty$, and observing that by construction the extensions of $\Q$ corresponding to the subgroups  $H$ of ${\rm im}(\rho)$ of index at most $N$  have this property.  Such a $q$ exists as a consequence of the theorems of Hermite-Minkowski and \v Cebotarev. Then by choice of $N,q,p$, using Theorem \ref{crucial}  and Corollary \ref{tobeused}, one sees that the projective image of  the image of a reduction of $\rho$ is either $\PSp_n(\F_{\ell^k})$ or $\PGSp_n(\F_{\ell^k})$ for some integer $k$.  By choosing $p$ appropriately we may ensure that $k$ is divisible by an integer $t$ chosen in advance.

The observations that if:
\begin{itemize}
\item a finite subgroup $G$ of $\GL_n({\overline \F_\ell})$ contains {\it deeply} embedded within 
it a  certain metacyclic subgroup,  then $G$ is forced to be {\it large}, and 
\item  the image of a global Galois representation can be made to contain such a metacyclic subgroup
by means of the Hermite-Minkowski theorem
\end{itemize}
we owe to  \cite{KW} in the case of $n=2$. Theorem \ref{crucial}  of this work generalizes the first
observation to all $n$.  The second observation can then be used in conjunction with 
automorphic methods to construct the required global Galois representations.

The main steps  to the proof of Theorem \ref{main} are:
\begin{enumerate}
\item A generalization of Lemma 6.3 of loc. cit. to any dimension (Theorem \ref{crucial}), and 
\item  Construction of self-dual, algebraic, regular cuspidal automorphic representations $\Pi$  on $\GL_n(\A_\Q)$, with $\A_\Q$ the adeles of $\Q$, with certain ramification properties: see Section \ref{last}.  The reader may consult  \cite{Cl} for the definition of regular and algebraic which is a condition on $\Pi_{\infty}$.
\end{enumerate}
Theorem 
\ref{crucial} might be of independent interest and be useful when extending the results of \cite{KW}. 

We indicate how we construct the $\Pi$'s: this  also allows us to introduce some necessary notation.

An expected source of  $\PSp_{n}(\bar \F_\ell)$-valued representations of $G_\Q$
are  self-dual automorphic representations $\Pi$ of $\GL_{2}(\A_\Q)$
which are regular algebraic at infinity and for which the exterior square $L$-function,
$L(s,\Lambda^2,\Pi)$, has a pole at $s=1$. 

For each place $v$ of $\Q$ we may attach to $\Pi_v$ its complex 
Langlands parameter $\sigma(\Pi_v)$  (we use the normalization of \cite{Cl}) which is a representation of the Weil-Deligne group $WD_v$ of $\Q_v$ with values in $\GL_{n}(\C)$. We may
regard this as valued in $\GL_{n}({\overline \Q_\ell})$ by choosing an isomorphism
$\C \simeq \overline \Q_\ell$. When $\Pi_v$ is unramified or supercuspidal,
$\sigma(\Pi_v)$ may be regarded as a representation of the Weil group $W_{\Q_v}
\subset WD_v$ of $\Q_v$; in fact, this will be the case at all finite places for the representations we construct.

The work in 
\cite{Kot}, \cite{Cl}, \cite{HT} attaches Galois representations to many such $\Pi$.
More precisely, if there is a finite  place $v$ such that $\Pi_v$ is discrete series, and $\Pi_\infty$ is regular and algebraic,  for every prime
$\ell$, there is an $\ell$-adic  Galois representation $\rho_{\Pi}:G_\Q \rightarrow \GL_{n}({\overline \Q_\ell})$ such that $\rho_{\Pi}|_{D_q}$ is isomorphic to  $\sigma(\Pi_v) \otimes |  \  \  |^{{1-n} \over 2}$  for all  primes  $q \neq p$ at which $\Pi_q$ is unramified or supercuspidal. 

We need to ensure certain ramification properties of $\Pi$ for this Galois representation to be of use to us. For this we give ourselves  the data of certain supercuspidal representations  $\pi_v$ of $\GL_{n}(\Q_v)$ for $v \in S$ a finite set of finite places  and a discrete series representation $\pi_\infty$ at  $ \infty$ with regular algebraic parameter. Then we have to construct a cuspidal automorphic representation $\Pi$ that is self-dual on $\GL_{n}(\A_{\Q})$, such that $\Pi$ is unramified outside $S$ and another place $w$ (which will typically be $\ell$), and $\Pi_v \simeq \pi_v$ for $v \in S \cup \{\infty\}$.

To construct self-dual automorphic representations on $\GL_n=\GL_{2m}$ which interpolate finitely many local self-dual representations directly seems rather subtle (see the remark 
on this point at the end of the paper). Instead we construct related 
generic cuspidal  automorphic representations on $\SO_{2m+1}(\A_\Q)$  
using Poincar\'e series (see Theorem \ref{global}) and transfer them to 
$\GL_{2m}(\A_\Q)$ using a known case of Langlands' principle of functoriality, namely  the forward lifting of Cogdell, Kim, Piatetski-Shapiro and Shahidi \cite{CKPSS}
that uses converse theorems. This accounts for the functoriality of the title (functoriality  is used  in some more of our references, e.g. \cite{Cl}).
The results of Jiang and Soudry \cite{JS1}, \cite{JS2} which  prove the local Langlands correspondence for generic supercuspidal representations of $\SO_{2m+1}(\Q_p)$, and that the lifts from $\SO_{2m+1}(\A_\Q)$ to 
$\GL_{2m}(\A_\Q)$ constructed in \cite{CKPSS} are functorial at all places, are crucial to us.

The $\ell$-adic representations $\rho_{\Pi}$ which arise this way from automorphic representations $\Pi$ on $\GL_n(\A_\Q)$ that are lifted from $\SO_{2m+1}(\A_\Q)$ come with a pairing $$ \rho \otimes \rho \rightarrow \Q_\ell(1-n).$$ It is expected, but probably not known in general, that this pairing can be chosen to be  symplectic. It also expected, but again not known in general, that if $\Pi$ is cuspidal, $\rho_{\Pi}$ is irreducible.  We use the fact that
the $\Pi$ we consider is such that $\sigma(\Pi_q)$ is an {\it irreducible} representation that preserves an  {\it alternating} form on $\overline \Q_\ell^n$, at some finite prime $q$,  to check this in the cases considered in this paper.

To summarize: we begin with a subgroup of $\Sp_{2m}(\bar\F_\ell)$ which can be realized as a Galois group over $\Q_q$ for a prime $q$ satisfying a certain condition of \v Cebotarev type.
We take the corresponding Weil group representation and use local Langlands for $\GL_n$ to construct a representation of $\GL_n(\Q_q)$.  We use inverse lifting 
to get a representation of $\SO_{2m+1}(\Q_q)$.  This becomes the factor at $q$ of an automorphic representation of $\SO_{2m+1}(\A_{\Q})$.  We then lift this to a self-dual representation on 
$\GL_{2m}(\A_{\Q})$, to which we associate a symplectic 
$\ell$-adic representation of $G_{\Q}$.  Thanks to known compatibilities, the restriction to $G_{\Q_q}$ of 
the reduction of this representation gives our original representation up to a twist.  Then a group theory argument (depending on the condition satisfied by $q$)
can be used to show that any subgroup of $\GSp_{2m}(\bar\F_\ell)$ 
which contains the image of
the specified image of $G_{\Q_q}$ is (up to conjugation and issues of center) of the form
$\Sp_{2m}(\F_{\ell^k})$
for some $k$ divisible by $t$.  

Some variant of this basic method might be made to work for other families of 
finite simple groups of Lie type.
It appears, however, that our poor control over which values of $k$ can be achieved is an 
unavoidable limitation of our technique, at least in its present form.  
We construct Galois representations by constructing cuspidal automorphic representations $\pi$ on $\GL_n(\A_\Q)$ using Poincar\'e series and the results of \cite{CKPSS}. Thus this  allows no control on the field of  definition of $\pi$. On the other hand by explicitly computing Hecke eigenvalues of cuspidal automorphic representations on $\SO_{2n+1}(\A_\Q)$, and choosing $\rho_q$ carefully, one  could in principle realize $\PSp_n(\F_{\ell^k})$ for specific values of $k$.

On the positive side, this method
does give good control of ramification.  In fact, all the Galois extensions of $\Q$ constructed in this
paper can be ramified only at $\ell$, $q$, and $\infty$.

We end our paper by proving that the $\ell$-adic Galois representations we construct, whose reductions mod $\ell$ enable us to prove Theorem \ref{main}, also have large images; namely their Zariski closure is $\GSp_n$.

We itemize the contents of the paper. In Section \ref{groups} we prove the 
group theoretic result Theorem \ref{crucial} that is key for us. In Section \ref{prelim} we fix the local Galois theoretic data that we need to realize as arising from a global Galois representation to prove Theorem \ref{main}.
 In Section \ref{Poincare} we prove Theorem \ref{global} which yields existence of generic cuspidal representations $\pi$ of a quasi-split group over $\Q$, with some  control on the  ramification of $\pi$, that interpolate finitely many   given local representations that are generic, integrable discrete series representations.
In Section \ref{proofs} we combine all the earlier work to prove Theorem \ref{main}. We end  with Section \ref{extra} that determines the Zariski closures of the images of the $\ell$-adic  Galois representations we construct.

We would like to thank Adri\'an Zenteno Guti\'errez for pointing out an error in the treatment of the characteristic $2$ case in the published version of this paper and 
Gabor Wiese and Luis Dieulefait for pointing out that the methods of this paper cannot distinguish between images of the form
$\PSp_n(\F_{\ell^k})$ and $\GSp_n(\F_{\ell^k})/\F_{\ell^k}^\times$.

\section{Some group theory}\label{groups}

Let $\Gamma$ be a group and $d\ge 2$ an integer.   We define $\Gamma^d$ as the intersection of
all normal subgroups of $\Gamma$ of index $\le d$.  

Let $n\ge 2$ be an integer and $p$ a prime congruent to $1$ (mod $n$).  By a
group \emph{of type $(n,p)$}, we mean any non-abelian homomorphic image
of any extension of $\Z/n\Z$ by $\Z/2p\Z$ such that $\Z/n\Z$ acts 
faithfully on $\Z/p\Z\subset \Z/2p\Z$.

These groups have the following property:

\begin{lemma}
If $G$ is a group of type $(n,p)$ and $\ell$ is a prime distinct from
$p$, then every faithful representation of $G$ over $\bar\F_\ell$
has dimension $\ge n$.
\end{lemma}

\begin{proof}
\label{G-np}
It suffices to prove that 
if $0\to\Z/2p\Z\to H\to \Z/n\Z\to 0$ and $\Z/n\Z$ acts faithfully
on $\Z/p\Z$, then every irreducible representation of $H$
has dimension $1$ or dimension $\ge n$.  The restriction of any 
such representation to $\Z/p\Z$ is a direct sum of characters
since $\ell\neq p$.  If every character is trivial,
then the original representation factors through an extension of $\Z/n\Z$ by 
$\Z/2\Z$, and such an extension is always abelian.  Otherwise, a non-trivial
character $\chi$ of $\Z/p\Z$ appears, so every character obtained by composing
$\chi$ with an automorphism of $\Z/p\Z$ coming from the action of $\Z/n\Z$ likewise appears.
As there are $n$ such characters, the original representation must have degree
$\ge n$.

\end{proof}

We can now state the theorem:

\begin{theorem}\label{crucial}
Let $n\ge 2$ be an integer.  There exist constants $N_d$ and $N_p$ such that if
$d > N_d$ is an integer, $p > N_p$ and $\ell$ are distinct primes, and 
$\Gamma\subset \GL_n(\bar\F_\ell)$ is a finite group such 
that $\Gamma^d$ contains a group of type $(n,p)$,
then there exists $g\in \GL_n(\bar\F_\ell)$
and $k\ge 1$ such that $g^{-1}\Gamma g$ is
one of the following:
\begin{enumerate}
\item A group containing $\SL_n(\F_{\ell^k})$ or $\SU_n(\F_{\ell^k})$
and contained in its normalizer.
\item A group containing $\Sp_n(\F_{\ell^k})$ and contained in its normalizer.
\end{enumerate}
\end{theorem}

\begin{proof}
By the main theorem of \cite{LP}, 
there exists a constant $J(n)$ depending only on $n$
such that every $\Gamma\subset\GL_n(\bar\F_\ell)$ has normal subgroups 
$\Gamma_1\subset\Gamma_2\subset \Gamma_3$ with the following properties:
\begin{itemize}
\item[(a)] $\Gamma_1$ is an $\ell$-group.
\item[(b)] $\Gamma_2/\Gamma_1$ is an abelian group of prime-to-$\ell$ order
\item[(c)] $\Gamma_3/\Gamma_2$ is isomorphic to a product 
$\Delta_1\times\cdots\times\Delta_r$ 
of finite simple groups of Lie type in characteristic
$\ell$.
\item[(d)] $\Gamma/\Gamma_3$ is of order $\le J(n)$.
\end{itemize}

If $\Gamma^d$ contains a subgroup of type $(n,p)$ for $d > J(n)$, 
then $\Gamma_3$ contains such a subgroup.
By Lemma~\ref{G-np}, any faithful $n$-dimensional representation
of a group of type $(n,p)$ in characteristic $\ell$ is irreducible
and tensor-indecomposable.  Thus the same is true for the representation
of $\Gamma_3$.  It follows that $\Gamma_1=\{1\}$.  We conclude that $\Gamma_3$ is
an abelian extension of $\Delta_1\times\cdots\times \Delta_r$.
This implies $r\ge 1$.

We have the following lemma:

\begin{lemma}
If $q\neq \ell$ is a prime, $\Delta$ is isomorphic to a product of
finite simple groups of Lie type in characteristic $\ell$,
and $\phi\colon \Delta\to\GL_n(\bar\F_q)$ is a homomorphism, 
then 
$$|\phi(\Delta)| \le \max(J(n), 25920)^{n/2}.$$
\end{lemma}

\begin{proof}
The image $\phi(\Delta)$ is again a product $\Delta_1\times\cdots\times \Delta_s$
of simple groups of Lie type in characteristic $\ell$, applying \cite{LP}
to $\phi(\Delta)$, and renumbering the $\Delta_i$ if necessary,
we may assume that there exists $t\le s$ so that
$$|\Delta_1|\cdots|\Delta_t| \le J(n)$$
and each $\Delta_i$ for $t<i\le s$ is of Lie type in characteristic $q$.  There are finitely many
finite simple groups which are of Lie type in two different characteristics, and the largest
is $U_4(\F_2)\cong \PSp_4(\F_3)$ \cite[p.~xv]{Atlas}, which is of order 25920.
Thus, 
$$|\phi(\Delta)| \le \max(J(n),25920)^s.$$
To bound $s$, we use the fact that every faithful irreducible representation of a product of
$k$ finite simple groups is an external tensor product of faithful representations of
these groups and therefore of degree $\ge 2^k$.  Every faithful representation
of $\Delta_1\times\cdots\times\Delta_s$ has, for each $i$ from $1$ to $s$, 
at least one irreducible factor which is faithful on $\Delta_i$.  Thus the dimension of
such a representation has degree at least $2^{k_1}+\cdots+2^{k_u}$ where
$k_1+\cdots+k_u = s$.  It follows that $2s\le n$.
\end{proof}

From this we can deduce the following:

\begin{lemma}
Let 
$$N_d = J(n) \max(J(n), 25920)^{n/2}$$
and $d>N_d$.
There exist normal subgroups $\{1\} = \Gamma'_1\subset\Gamma'_2\subset\Gamma'_3$
of $\Gamma$ satisfying conditions (a)--(c) above together with two additional conditions:
$\Gamma'_2$ lies in the center of $\Gamma'_3$ and 
$\Gamma^d\subset \Gamma'_3$.
\end{lemma}

\begin{proof}
Without loss of generality we may assume that 
$\Gamma_3/\Gamma_2$ is non-trivial,
in which case it is a non-trivial product of groups of Lie type in characteristic $\ell$.
Let $q$ denote a prime dividing the order of $\Gamma_2$.
Thus $q\neq\ell$.  Let $\Gamma_2[q]$ and $\Gamma_2[q^\infty]$ denote the
kernel of multiplication by $q$ and the $q$-Sylow subgroup respectively.
As $\Gamma_2[q]$ is an elementary abelian $q$-group contained in
$\GL_n(\bar\F_\ell)$, its dimension as $\F_q$-vector space is $\le n$.
By the preceding lemma, the image of the homomorphism
$$\phi\colon \Gamma_3/\Gamma_2\to \Aut\Gamma_2[q]\subset \GL_n(\F_{q})$$
giving the action of $\Gamma_3/\Gamma_2$ on $\Gamma_2[q]$ has order bounded
above by
$\max(J(n), 25920)^{n/2}$.  Let $\Gamma_{3,q}$ denote the preimage in $\Gamma_3$ of
$\ker\phi$.  As $\Gamma_2[q]$ is normal in $\Gamma$, we see that $\Gamma_{3,q}$ is a normal subgroup of $\Gamma$ of index $\le J(n)|\ker\phi|<d$.
Let $\Gamma'_3$ denote the intersection of $\Gamma_{3,q}$ over all primes $q$ dividing the order
of $\Gamma_2$.  Then $\Gamma^d\subset\Gamma'_3$, and
$\Gamma'_3/\Gamma_2$ is a normal subgroup of a product of groups of Lie type in
characteristic $\ell$ and is therefore again such a product.
Its action on
each $\Gamma_2[q^\infty]$ is trivial since
$\ker \Aut\Gamma_2[q^\infty]\to \Aut\Gamma_2[q]$ is a $q$-group.
Therefore, its action on $\Gamma_2$ is trivial.  Setting $\Gamma'_2=\Gamma_2$, we get the lemma.
\end{proof}

Redefining $\Gamma_i := \Gamma'_i$, we may assume that (a)--(c) hold together with the
condition $\Gamma^d\subset \Gamma_3$, and we proceed on the hypothesis that
$\Gamma^d$ contains a subgroup of type $(n,p)$. 
A consequence of the preceding lemma is that the representation of $\Gamma_3$
can be regarded as a projective representation of 
$\Delta_1\times\cdots\times \Delta_r$, and as such, it is 
tensor-indecomposable.  It follows that $r=1$.
As $\Delta_1$ is a group of Lie type in characteristic $\ell$, there
exists a simply connected almost simple 
algebraic group $D/\bar\F_\ell$ and a Frobenius map
$F\colon D\to D$ such that $\Delta_1$ is isomorphic to the quotient of
$D(\bar\F_\ell)^F$ by its center.  Moreover, $D(\bar\F_\ell)^F$
is the universal central extension of $\Delta_1$, so the projective
representation $\Delta_1\to\PGL_n(\bar\F_\ell)$ lifts to an irreducible linear representation
$D(\bar\F_\ell)^F\to\GL_n(\bar\F_\ell)$.
By a well-known theorem of Steinberg
\cite[13.1]{Steinberg}, the irreducible representations of
$D(\bar\F_\ell)^F$ over $\bar\F_\ell$ extend to irreducible representations of the algebraic group
$D$.  Thus we have a non-trivial representation
$\rho\colon D\to \GL_n$.  In particular, $\dim D\le n^2$. 

As $\Gamma_3$ contains a group of type $(n,p)$
it contains an element $x$ of order $p$ and an element $y$
such that $y^{-1} x y = x^a$, where $a\in \Z/p\Z$ is an element of order $n$.
Thus $x$ and $y^{-1} x y$ are commuting elements of order $p$.
Let $\tilde x,\tilde y\in D(\bar\F_\ell)$ lie over $x,y\in\Delta_1$.
Let $x_1 = \tilde x$ and $x_2=\tilde y^{-1}\tilde x\tilde y$.  Thus
$x_2$ lies over $x^a$.  It follows that the commutator of $x_1$ and $x_2$ lies in the center of $D$.

We can therefore make use of the following lemma:

\begin{lemma}
\label{torify}
Let $G$ be an semisimple algebraic group 
over an algebraically closed field $F$.  
Then there exists a constant $N$ depending only on $\dim G$ such that if 
$p>N$ is prime and $p\neq 0$ in $F$, then for any two elements $x_1,x_2\in G(F)$ of order $p$ 
whose commutator lies in the center of $G$ 
there exists a maximal torus $T$ such that $x_1,x_2\in T(F)$.
\end{lemma}

\begin{proof}
We use induction on $\dim G$.
If $x_1$ and $x_2$ lie in the center of $G$, any maximal torus will do.  
Without loss of generality we assume
that $x_1$ is not central.
If $\tilde x_i$ denotes a preimage of $x_i$ in $\tilde G(F)$, where $\tilde G$ is the universal cover of $G$,
then $\tilde x_1\tilde x_2 = z(x_1,x_2)\tilde x_2\tilde x_1$, where $z(x_1,x_2)$ lies in 
the center of $\tilde G$.  If $p$ is greater than the order of the center of $\tilde G$,
this implies that $\tilde x_1$ and $\tilde x_2$ commute, so $x_1$ and $x_2$ lie in the image $H$ in $G$
of the centralizer $Z_{\tilde G}(\tilde x_1)$.  
Note that $x_1$ is semisimple due to its order, so $\tilde x_1$ is semisimple, and
by Steinberg's theorem,
$Z_{\tilde G}(\tilde x_1)$ is a connected reductive group.  As $H$ is connected
and reductive, it can be written $H = H' Z$, where the derived group $H'$ of $H$ is semisimple
and $Z$ is the identity component of the center of $H$, which is a torus.  
Let $x_i=x'_iz_i$ for $i=1,2$.
By the induction hypothesis, $x'_1$ and $x'_2$
lie in a common maximal torus $T'$ of $H'$, and setting $T = T' Z$, the lemma follows by induction.
\end{proof}
 
It follows that there exists a maximal torus $T$ in $D$ such that
$x_1$ and $x_2$ both lie in $T(\bar\F_\ell)$.  
By a well-known theorem \cite[\S3.1]{Humphreys},
there exists $w$ in the normalizer of $T$ such that 
$$w^{-1} x_1 w = x_2 = x_1^a z.$$
As $x_1$ and $x$ have the same image in $\PGL_n(\bar\F_\ell)$,
$$\rho(x_1) \sim \omega
\begin{pmatrix} 
\lambda&0&\cdots&0\\ 0&\lambda^a&\cdots&0\\ \vdots&\vdots&\ddots&0\\0&0&\cdots&\lambda^{a^{n-1}}
\end{pmatrix}
$$
for some $\omega$.  This implies that the characters of $\rho$ with respect to
$T$ are pairwise distinct.

Conjugation by $w$ permutes the weights of $\rho$ cyclically.
In particular, the Weyl group acts transitively on the weights,
so $\rho$ is miniscule.  By the classification of miniscule representations
\cite[Annexe]{Serre}, one of the following must hold:

\begin{enumerate}
\item $D=\SL_m$ and $\rho$ is a fundamental representation.
\item $D=\Sp_n$, and $\rho$ is the natural representation.
\item $D=\Spin_{2m}$, $n=2m$, and $\rho$ is the natural representation
of $\SO_n$.
\item $D=\Spin_{2m}$, $n=2^{m-1}$, and $\rho$ is a semispin representation.
\item $D=\Spin_{2m+1}$, $n=2^m$, and $\rho$ is the spin representation.
\item $D=E_6$ and $n=27$.
\item $D=E_7$ and $n=56$.
\end{enumerate}

In case (1), $\rho$ must be the natural representation or its dual because 
no permutation of an $m$-element set $S$ generates a group acting
transitively on the set of $k$-element subsets of $S$ when $2\le k\le m-2$.
In case (3), the image of $w$ in $S_m$ generates a group which acts transitively on $\{1,2,\ldots,m\}$.
However, all elements in the Weyl group of the root system of $D_m$ mapping to an $m$-cycle in $S_m$
are conjugate to one another, and all are therefore of order $m$.  For the remaining spin cases, 
we can take an odd power $w^k$ of $w$ whose order is of the form $2^r$ and such that $w^k$ acts transitively on the $2^m$ or $2^{m-1}$ weights
of a spin or semispin representation.  The image of $w^k$ in $S_m$ has order $2^r$ or $2^{r-1}$, so $2^{r-1} \le m$.  Thus, $2^{m-2}\le m$, which means $m\le 4$.
An examination of cases shows that this can occur for the spin representation of $\Spin_5$ and either semispin representation of $\Spin_6$.  The resulting cases
duplicate (2) for $n=4$ and (1) for $m=4$ respectively.
Cases (6) and (7) cannot occur because the Weyl groups of $E_6$ and $E_7$ have no elements whose order is divisible by $27$ or $56$ respectively \cite{Atlas}.

We conclude that it suffices to consider the cases
\begin{enumerate}
\item $D=\SL_n$ and $\rho$ is the natural representation.
\item $D=\Sp_n$ and $\rho$ is the natural representation.
\end{enumerate}
In case (1), $D(\bar\F_\ell)^F$ is of the form $\SL_n(\F_{\ell^k})$
or $\SU_n(\F_{\ell^k})$.
In case (2), $D(\bar\F_\ell)^F=\Sp_n(\F_{\ell^k})$.

Now,
$$[\Gamma_3,\Gamma_3]=[\rho(D(\bar\F_\ell)^F),\rho(D(\bar\F_\ell)^F)]
=\rho(D(\bar\F_\ell)^F).$$
The possibilities for $\rho(D(\bar\F_\ell)^F)$ are
$\SL_n(\F_{\ell^k})$, $\SU_n(\F_{\ell^k})$, and $\Sp_n(\F_{\ell^k})$.
As
$$[\Gamma_3,\Gamma_3]\subset\Gamma\subset 
\mathrm{Norm}_{\GL_n(\bar\F_\ell)}(\Gamma_3),$$
we have the theorem.
\end{proof}

\begin{cor}\label{tobeused}
Under the hypotheses of Theorem~\ref{crucial}, if $\Gamma\subset\GSp_n(\bar\F_\ell)$,
and $\bar\Gamma$ denotes the image of $\Gamma$ in $\PGL_n(\bar\F_\ell)$, then there exists
$\bar g\in \PGL_n(\bar\F_\ell)$ and a positive integer $k$ such that 
$$\bar g^{-1}\bar \Gamma\bar g\in\{\PSp_n(\F_{\ell^k}),\GSp_n(\F_{\ell^k})/\F_{\ell^k}^\times\}.$$
\end{cor}

\begin{proof}
If $g^{-1}\Gamma g$ contains $\SL_n(\F_{\ell^k})$ or $\SU_n(\F_{\ell^k})$,
then one of these groups has an $n$-dimensional
symplectic representation.  When $n=2$, $\SL_n$, $\SU_n$, and $\Sp_n$ all coincide, so $g^{-1}\Gamma g$ contains $\Sp_n(\F_{\ell^k})$.  For $n\ge 3$,
from Steinberg's theorem, it follows that the algebraic group $\SL_n$
has a non-trivial self-dual $n$-dimensional representation defined over $\bar\F_\ell$
which maps the fixed points of a Frobenius
map into $\Sp_n(\bar\F_\ell)$.  Of course $\SL_n$ has no non-trivial self-dual representation of dimension $n$ when $n\ge 3$.  

In any case, by Theorem~\ref{crucial}, $g^{-1}\Gamma g$ is trapped between
$\Sp_n(\F_{\ell^k})$ and its normalizer in $\GL_n(\bar\F_\ell)$.   To compute the normalizer,
we first note that $\Sp_n(\F_{\ell^k})$ has no non-trivial graph automorphisms, so its
outer automorphism group is the semidirect product of the group of diagonal automorphisms
$\Z/2\Z$ (or $\{0\}$ if $\ell=2$) and the group of field automorphisms $\Z/k\Z$.  

Non-trivial field automorphisms
never preserve the character of the $n$-dimensional representation of $\Sp_n(\F_{\ell^k})$.
For $\ell^k\neq 4$, we can see this by noting that, by a counting argument, $\F_{\ell^k}$ contains an element $\alpha$
such that $\alpha+\alpha^{-1}$ is not contained in any proper subfield, and there exists
an element of $\Sp_n(\F_{\ell^k})$ with eigenvalues $1,1,\ldots,1,\alpha,\alpha^{-1}$.
For $n\ge 4$, there exists an element $\alpha$, of $\bar\F_4$ of order $17$ and an element
of $\Sp_n(\F_4)$ with eigenvalues $1,1,\ldots,1,\alpha,\alpha^4,\alpha^{-4},\alpha^{-1}$ and therefore with trace in $\F_4\setminus\F_2$.  Finally, $\SL_2(\F_4)$ contains the element
$\binom{1\,\omega^{\phantom2}}{1\,\omega^2}$ with trace $\omega\not\in\F_2$.  Thus,
$$[N_{\GL_n(\bar\F_\ell)}\Sp_n(\F_{\ell^k}):\Sp_n(\F_{\ell^k})\bar\F_\ell^\times]
\le\begin{cases}2&\text{if $\ell$ is odd,} \\ 1&\text{if $\ell=2$.} \\ \end{cases}$$

On the other hand, when $\ell$ is odd, 
$\GSp_n(\F_{\ell^k})\subset \Sp_n(\F_{\ell^k})$ contains elements which act on
$\Sp_n(\F_{\ell^k})$ by the non-trivial diagonal automorphism.  By Schur's lemma, 
the normalizer of $\Sp_n(\F_{\ell^k})$ in $\GL_n(\bar\F_\ell)$ is therefore
$\GSp_n(\F_{\ell^k})\bar\F_\ell^\times$.  

\end{proof}

\noindent{\bf Remark:} We indicate  how the proof of  Theorem \ref{crucial} is related to that of  Lemma
6.3 of \cite{KW}. There it proved that every subgroup $G$ of $\PGL_2(\bar \F_\ell)$, whose projective image has the property that every index 2 subgroup of $G$ contains the dihedral group of order $2p$ with $p>5$ 
a prime, is conjugate to a subgroup that is trapped between 
$\PSL_2(\F_{\ell^k})$ and $\PGL_2(\F_{\ell^k})$ for some integer $k$. This is proved using Dickson's theorem. The role of Dickson's theorem here is played by the results of \cite{LP}.

\section{A few preliminaries for the proof of Theorem \ref{main}}\label{prelim}

\subsection{A tamely ramified symplectic local parameter at $q$ of dimension $n$}

Let $p,q>n$ be distinct odd primes, such that the order of $q$ mod $p$ is $n=2m$.
Consider the degree $n$ unramified extension $\Q_{q^n}$ of $\Q_q$.
We consider a character $\chi:\Q_{q^n}^\times \simeq \mu_{q^n-1} \times U_1 \times q^{\Z} \rightarrow \aQl^\times$ such that 
\begin{itemize} \item the order of $\chi$ is $2p$ \item $\chi|_{\mu_{q^n-1} \times U_1}$ is of order $p$ \item $\chi(q)=-1$.\end{itemize}  We call such a $\chi$ a tame symplectic character of $\Q_q$ of degree $n$ and order $2p$.  By local class field theory, we can regard $\chi$ as a character of
$G_{\Q_{q^n}}$.  (We normalize the isomorphism of class field theory by sending 
a uniformizer to an arithmetic Frobenius.)

Consider $\rho_q:G_{\Q_q} \rightarrow \GL_n(\aQl)$ that is given by 
${\rm Ind}_{\Q_{q^n}}^{\Q_q}\chi$.

The following is easily deduced from Theorem 1 of \cite{Moy}:

\begin{prop}\label{moy}
The representation $\rho_q$ is irreducible and 
symplectic, and thus it can be conjugated to take values in $\Sp_n(\aQl)$. 
\end{prop}

\begin{proof}
 The irreducibility follows from the fact that the order of $\chi$ is $2p$ and the order of $q$ mod $2p$ is $n$. This ensures that the characters $\chi,\chi^q,\chi^{q^2},\cdots,\chi^{q^{n-1}}$ are all distinct. Also note that $\chi|_{\Q_{q^m}^\times}$ is unramified (i.e., trivial on the units of $\Q_{q^m}$) and
 of order 2. Then Theorem 1 of \cite{Moy} proves that $\rho_q$ is symplectic.
\end{proof}

We assume  $p \neq \ell$.  The image of the reduction of  $\rho_q$ in $\GL_n({\overline \F_\ell})$ is a group of type $(n,p)$.  It acts irreducibly on  $\bar\F_\ell^n$ and preserves up to scalars a unique bilinear form which is
necessarily non-degenerate and alternating.

\subsection{Some lemmas}

Next we recall some well-known facts concerning the values of
cyclotomic polynomials.  Let $R_n$ denote the set of primitive complex
$n$th roots of unity, and
$$\Phi_n(x) = \prod_{\zeta\in R_n} (x-\zeta).$$
If $a$ is an integer, $n$ a positive integer, and $p$ a prime dividing $\Phi_n(a)$,
then either the class of $a$ in $\F_p^\times$
has order exactly $n$ or $p$ divides $n$ \cite[Lemma 2.9]{Wa}.
In the former case, $p$ cannot divide $\Phi_d(a)$ for any proper divisor $d$ of $n$.
In the latter case, we have the following result:

\begin{lemma}
If $n\ge 3$, $a\in\Z$, and $p$ is a prime dividing $n$, then $p^2$ does not
divide $\Phi_n(a)$.
\end{lemma}

\begin{proof}
Suppose first that $p=2$ and $n = 2^k$ for $k\ge 2$ an integer.  Then
$$\Phi_n(a) = a^{2^{k-1}}+1 = (a^{2^{k-2}})^2+1 \not\equiv 0 \pmod4.$$
If $p=2$ and $n$ has an odd prime divisor $q$, then $\Phi_n(x)$ divides
$\Phi_q(x^{n/q})$ in $\Z[x]$, so $\Phi_n(a)$ divides
$$1+a^{\frac nq}+a^{\frac {2n}{q}}+\cdots+a^{\frac{(q-1)n}{q}}\equiv 1\pmod2.$$
Finally, if $p$ is odd, $\Phi_n(a)$ divides $\Phi_p(a^{n/p})$  As $\Phi_p(x+1)$ is an Eisenstein
polynomial, evaluating $\Phi_p$ at an integer cannot give a multiple of $p^2$.
\end{proof}

From this we easily deduce the following:

\begin{lemma}
\label{cyclotomic-primes}
If $a\ge 3$ and $n\ge 3$ or $a=2$ and $n\ge 7$, then $\Phi_n(a)$ has a prime divisor $q$ such that 
the class of $a$ in $\F_q^\times$ has order exactly $n$.
\end{lemma}

\begin{proof}
It suffices to prove that $|\Phi_n(a)| > n$.   We first consider the case $a\ge 3$.
Then $|a-\zeta|> 2$ for every $\zeta\in R_n$, so we have
$|\Phi_n(a)| >  2^{\phi(n)}$.   For every prime power $P$ except for $P=2$, we have $\phi(P)\ge\sqrt P$.  As $\phi$ is multiplicative, for all $n\ge 1$, we have $\phi(n)\ge \sqrt{n/2}$.  For $x>2$,
$\log_2(x) < \sqrt{x/2}$, so $2^{\phi(n)} > n$ for all $n\ge 3$.

For $a=2$, we write
$$\log \Phi_n(x) = \sum_{d\mid n} \mu\bigl(\frac nd\bigr)\log(x^d-1).$$
As
$$|\log (2^d-1) - d\log 2| \le 2^{-d}(1-2^{-d})^{-1} \le 2^{1-d},$$
we have
$$|\log \Phi_n(2) - \phi(n)\log 2|\le \sum_{d=1}^\infty 2^{1-d} = 2.$$
For $n\ge 181$, we have
$$\phi(n)-2\ge \sqrt{n/2} - 2 > \log_2(n),$$
so we need only check $n\le 180$.  The only values $n\ge 7$ for which 
$\phi(n)-2\le \log_2(n)$ are $n=8,10,12,18$ for which $\Phi_n(2)$ has prime divisor
$17,11,13,19$ respectively.
\end{proof}

Now we can construct the primes $p$ and $q$ needed for the main theorem.

\begin{lemma}
\label{pq}
Given an even integer $n\ge 2$, a prime $\ell$, a finite Galois extension $K/\Q$, 
and positive integers $t$ and $N$, there exist primes 
$p$ and $q$ with the following properties:
\begin{enumerate}
\item The primes $p$, $q$, and $\ell$ are all distinct.
\item The prime $p$ is greater than $N$.
\item The prime $q$ splits in $K$.
\item The order of the image of $q$ in $\F_p^\times$ is exactly $n$.
\item If $\F$ is any finite field in characteristic $\ell$ and $\GSp_n(\F)$ contains an element of
order $p$, then $\F$ contains $\F_{\ell^t}$.
\end{enumerate}
\end{lemma}

\begin{proof}
Let $n=2m$.  Let $u>0$ denote a multiple of $t\cdot (m-1)!$.
Let $p$ be chosen to be a prime dividing $\Phi_{nu}(\ell)$ and therefore
$\Phi_n(\ell^u)$.  
We can make $p$ as large as we please by choosing $u$ sufficiently large.
We may therefore assume that $p>\max(n,\ell,N)$ and $K/\Q$ is not ramified at $p$.
Let $q\neq \ell$ be a prime congruent to $\ell^u$ (mod $p$), split in $K$,
and greater than $n$.  As $\Q(\zeta_p)$ and $K$ are linearly disjoint over $\Q$,
the \v Cebotarev density theorem guarantees the existence of such a prime.
As $p\neq q$, the first property is satisfied.  
The second and third properties are built into the definitions of $p$ and $q$ respectively.  
As $p$ does not divide $n$ and
$$\Phi_n(q)\equiv \Phi_n(\ell^u)\equiv 0\pmod p,$$
the fourth property is satisfied.  For the last property, we note that
$$|\GSp_n(\F_{\ell^k})| = (\ell^k-1)\ell^{km^2}\prod_{i=1}^m (\ell^{2ik}-1).$$
If $\GSp_n(\F_{\ell^k})$ has an element of order $p$,
then $p$ divides $\ell^{2ik}-1$ for some $i$ between $1$ and $m$, which means that
 the order of $\ell$ in $\F_p^\times$ divides $2k i$ for some $i\le m$
and therefore divides $2k\cdot m!$.  We know that the order is in fact $nu$, which is an integral
multiple of $2t\cdot m!$, so $t$ divides $k$, as claimed.

\end{proof}

\subsection{Fixing Galois-theoretic  data}

Let $t$ be a given positive integer.  We may freely replace $t$ by any positive multiple, so
without loss of generality we assume that $t$ is divisible by $n$.
We define $N$ to be the $N_d$ of Theorem \ref{crucial}, and let $K$ denote the compositum of
all Galois extensions of $\Q$ of degree $\le N$ which are ramified only over $\ell$
and $\infty$.  By the Hermite-Minkowski theorem, $K$ is a number field.
We define $p$ and $q$ via Lemma~\ref{pq}
and consider
the representation $\rho_{q}={\rm Ind}_{\Q_{q^n}}^{\Q_q}\chi:G_{\Q_{q}} \rightarrow \Sp_n(\aQl)$ for $\chi$ a tame, symplectic character of 
$\Q_{q}$ of degree $n$ and order $2p$.
Note that $\chi(I_q)$ has order $p$.

\section{Globalizing discrete series}\label{Poincare}

In this section we show how to construct a global, generic cuspidal representation with 
desired local components. The precise result is contained in Theorem \ref{global}. 

\subsection{Poincar\' e Series}

Let $G$ be a quasi-split simple algebraic 
  group over $\mathbb Q$ corresponding to a (simple)
root system. The group 
 $K_{p}=G(\mathbb Z_{p})$ 
is a hyperspecial maximal compact subgroup in $G(\mathbb Q_{p})$ for almost all 
primes. We assume that the Lie group $G(\mathbb R)$ has discrete series representations. This 
condition determines the quasi-split $G(\mathbb R)$, up to an isogeny. 

Let $(\pi,H)$ be an \underline{integrable} discrete series of $G(\mathbb R)$ on a Hilbert space $H$.  
 Fix $K$, a maximal compact subgroup in $G(\mathbb R)$. The space $H_{K}$ of $K$-finite 
vectors in $H$ is an irreducible $(\mathfrak g,K)$-module. 
Let $f=f_{\infty}\otimes_{p} f_{p}$ be a function on $G(\mathbb A)$ 
such that $f_{p}$ is compactly supported for every prime $p$ and $f_{p}$ is 
equal to the characteristic function of $K_{p}$ for almost all primes. Moreover,
 $f_{\infty}$ is a matrix coefficient of the integrable discrete series. More 
 precisely, let $\langle v,w\rangle$ denote the inner product on $H$. 
 For our purposes the matrix coefficient $f_{\infty}$
 is a function 
 $$ 
 f_{\infty}(g)=\langle \pi(g)w, \pi(g_{1})v\rangle
 $$ 
 where $v$ and $w$ are $K$-finite vectors in $H$ and $g_{1}$ is an element in 
 $G(\mathbb R)$.  Let $Z(\mathfrak g)$ be the center of the enveloping 
 algebra of $\mathfrak g$. The function $f$ satisfies: 
 \begin{itemize}
 \item $f$ is in $L^{1}(G(\mathbb A))$. 
 \item $f$ is right $K$-finite. 
 \item $f$ is an eigenfunction of $Z(\mathfrak g)$.
 \end{itemize}
 It follows (see \cite{Bo2}) that the Poincar\' e series 
 $$ 
 P_{f}(g)=\sum_{\gamma\in G(\mathbb Q)} f(\gamma g)
 $$ 
converges absolutely and in the $C^{\infty}$-topology 
 to a smooth function on $G(\mathbb Q)\backslash G(\mathbb A)$. 
 In particular, the series converges uniformly on compact sets in $G(\mathbb A)$. 
 This function is cuspidal. That is, for every parabolic subgroup $P=MN$ defined over 
 $\mathbb Q$, the constant term
 $$ 
c_{N}(P_{f})(g)=\int_{N(\mathbb Q)\backslash N(\mathbb A)} P_{f}(ng) ~dn
$$
vanishes. This is easy to verify. Indeed, the
   matrix coefficient $f_{\infty}$ lies in the sub-space of cusp forms
 $\mathcal C_{0}(G(\mathbb R))$ in the Schwarz space 
 $\mathcal C(\mathbb R)$ of Harish-Chandra. In particular, 
 $$ 
 \int_{N(\mathbb R)} f_{\infty}(ng) ~dn=0.
 $$
 Since the Poincar\' e series is uniformly convergent on compact sets and the integral 
 defining the constant term is taken over a compact set, we can switch the order 
 of integration and summation to obtain 
$$
c_{N}(P_{f})(g)=
\sum_{\gamma\in G(\mathbb Q)/N(\mathbb Q)}
 \prod_{v} \int_{N_{v}}f_{v}(\gamma ng_{v})~dn,
$$
where we have abbreviated $N_{v}=N(\mathbb Q_{p})$ if $v=p$ and $N_{v}=N(\mathbb R)$ if $v=\infty$. 
It follows that $c_{N}(P_{f})=0$ since the local integral vanishes for $v=\infty$. 
 
 For every $X$ in the Lie algebra $\mathfrak g$ 
let $R_{X}$ denote the natural right action of $X$ on smooth functions on 
$G(\mathbb R)$. Then 
$$ 
R_{X}(P_{f})=P_{Xf}
$$
where $Xf(g)=\langle \pi(g)\pi(X)w,\pi(g_{1})v\rangle$. It follows that, 
by fixing the $K$-finite $v$ in $H$, the map 
$$ 
w\mapsto P_{f}
$$
is an intertwining map---in the sense of $(\mathfrak g,K)$-modules---from $H_{K}$ into 
$C^{\infty}_{0}(G(\mathbb Q)\backslash G(\mathbb A))_{K}$. In addition, 
for any prime $q$, the local factor
$f_{q}$ can be taken to be a matrix coefficient of a supercuspidal representation $\pi_{q}$
of $G(\mathbb Q_{q})$. Then the Poincar\' e series, if non-vanishing,  will give rise to a 
cuspidal automorphic representation which has the integrable discrete 
series at the real place and the supercuspidal representation $\pi_{q}$ as a local 
factor at the prime $q$.

 \subsection{Genericity of Poincar\' e series}
 
  Let $N$ be the unipotent radical of a Borel subgroup $B$ of $G$, defined over 
 $\mathbb Q$. Fix $\psi$ a Whittaker character of 
 $N(\mathbb A)$ trivial on $N(\mathbb Q)$.  Recall that an automorphic 
 representation $\pi$-generic if 
$$ 
W_{\psi}(\phi)=\int_{N(\mathbb Q)\backslash N(\mathbb A)} 
\phi(n)\psi(n) ~dn \neq 0
$$ 
for some (smooth) $\phi$ in $\pi$.

 Fix two finite and disjoint sets of places: $D$, containing $\infty$ 
and perhaps nothing else,  and
$S$, a non-empty set of primes such that $G$ is unramified at all primes $p$ 
not in $D\cup S$. This means that 
  $K_{p}=G(\mathbb Z_{p})$ is a 
 hyperspecial maximal compact subgroup of $G(\mathbb Q_{p})$. 
 For $G$ split, $S$ could consist of 
 only one prime. 
 %The primes in $S$ are called \emph{sacrificial primes}. 
 %Comment: we don't use this terminology elsewhere in the paper.
 In this section we shall show how Poincar\' e series gives a globally $\psi$-generic 
 (and thus non-zero) cuspidal automorphic representation $\pi$ such that
 
 \begin{itemize}
 \item $\pi_{\infty}$ is a 
 (given) generic integrable discrete series representation. 
 \item $\pi_{q}$ is a (given) generic supercuspidal representation for every prime $q$ in $D$. 
 \item $\pi_{p}$ is unramified for all $p$ not in $S\cup D$. 
 \end{itemize}
 
  We assume, as we can, that 
 $\psi$ is trivial on $N_{p}\cap K_{p}$ for every prime $p$ not in $S$.

  The Poincar\' e series is constructed as follows: Let 
  $f=\otimes_{v} f_{v} $ be a function on $G(\mathbb A)$
  such that: 
  
  \begin{itemize}
  \item $f_{\infty}$ is a matrix coefficient of the generic
   integrable discrete series $\pi_{\infty}$. 
   \item $f_{q}$ is a (compactly supported) matrix coefficient of the generic 
   supercuspidal representation $\pi_{q}$ for every prime $q$ in $D$. 
  \item $f_{p}$ is the characteristic function of $K_{p}=G(\mathbb Z_{p})$
   for all $p$ not in $S\cup D$. 
  \end{itemize} 
 
 We shall specify the local components $f_{\ell}$ for $\ell$ in $S$ in a moment.
 The idea is to show that for some choice of $f_{\ell}$, the Poincare series 
 is generic.  Let $B^{-}$ be the Borel subgroup opposite to $B$. 
 For every prime $\ell$ in $S$ pick a decreasing sequence $K_{\ell}^{m}$ 
 of open compact subgroups $K_{\ell}$ such that 
 \begin{itemize}
 \item $K^{m}_{\ell}\cap N_{\ell}$ is independent of $m$ and 
 $\psi$ is trivial on it. 
 \item  $K^{m}_{\ell}$ admits a parahoric factorization 
 $$ 
 K^{m}_{\ell}=(K^{m}_{\ell}\cap B^{-}_{\ell})\dot 
 (K^{m}_{\ell}\cap N_{\ell}). 
 $$
 \item $\lim_{m \rightarrow \infty}K^{m}_{\ell}\cap B_{\ell}^{-}=1$
 meaning that 
 $$ 
 \cap_{m=1}^{\infty}(K^{m}_{\ell}\cap B^{-}_{\ell})=\{1\}.
 $$
 \end{itemize}
 
Let $f^{m}$ be the 
function on $G(\mathbb A)$ which has the local factors outside $S$ independent of 
$m$ and as specified above, and 
$f_{\ell}^{m}$ the characteristic function of $K_{\ell}^{m}$ for all $\ell$ in $S$.
We shall show that $W_{\psi}(P_{f^{m}})(g)\neq 0$ for a sufficiently large $m$. 
In fact we can accomplish this with $g$ in $G(\mathbb A)$ such that $g_{p}=1$ 
for all $p$ not in $D$. 
Now to the proof. In order to save notation, assume that $D=\{\infty\}$ and 
$S=\{\ell\}$. 
Since the Poincar\' e series $P_{f^{m}}$ is uniformly convergent
on compact sets, we can switch the order of integration and summation. Assuming also that 
$g_{p}=1$ for all finite primes we have (an absolutely convergent series)
$$
W_{\psi}(P_{f^{m}})(g)=\sum_{\gamma\in G(\mathbb Q)/N(\mathbb Q)}
\int_{N_{\infty}} f_{\infty}(\gamma ng_{\infty})\psi(n)~dn  
 \prod_{p} \int_{N_{p}}f^{m}_{p}(\gamma n)\psi(n)~dn.
$$
Let $\Phi(f^{m}_{p},\gamma)$ denote the local factor in the above product. For a given 
$\gamma$ in $G(\mathbb Q)$, as $m$ varies, only the factor at $p=\ell$ could 
possibly change.

\begin{lemma}  Let $\gamma$ in $G(\mathbb Q)$ such that 
$\Phi(f_{\ell}^{m},\gamma)\neq 0$. Then
$$ 
\Phi(f_{\ell}^{m},\gamma)=\Phi(f_{\ell}^{1},\gamma).
$$
\end{lemma}
\begin{proof} If $\Phi(f_{\ell}^{m},\gamma)\neq 0$ then for some $n$ in $N_{\ell}$, 
$\gamma n\in K_{\ell}^{m}$ which means that $\gamma$ can be written as 
 $$
 \gamma=k_{\gamma}n_{\gamma}
 $$
 for some $k_{\gamma}$ in $K_{\ell}^{m}\cap B^{-}_{\ell}$ and $n_{\gamma}$
 in $N_{\ell}$. A trivial computation now shows that 
 $$ 
 \Phi(f_{\ell}^{1},\gamma)
 =vol(K_{\ell}^{1}\cap N_{\ell}) \cdot \psi(n_{\gamma})^{-1}
 =vol(K_{\ell}^{m}\cap N_{\ell}) \cdot \psi(n_{\gamma})^{-1}
 =\Phi(f_{\ell}^{m},\gamma). 
 $$
 \end{proof}
 
 The above lemma shows that the terms in the series
 $W_{\psi}(P_{f^{m}})(g)$ ($g$ is here fixed and trivial at all finite places) are 
 the same as the terms in the series $W_{\psi}(P_{f^{1}})(g)$ except we take 
 only $\gamma$ contained in 
 $$
  (K_{\ell}^{m}\cap B^{-}_{\ell})\cdot N_{\ell}.
 $$
 As $m$ goes to infinity, we are reduced to $\gamma$ which sit in $N_{\ell}$, 
 that is, in $N(\mathbb Q)$.  Since $\gamma$ a coset in $G(\mathbb Q)/N(\mathbb Q)$, 
 we can take $\gamma=1$ and the limit is equal to  
$$
\lim_{m\rightarrow \infty}W_{\psi}(P_{f^{m}})(g)= 
vol(K_{\ell}^{1}\cap N_{\ell})
\int_{N(\mathbb R)} f_{\infty}(ng_{\infty})\psi(n)~dn.
$$
(Here we assumed, as we could, that $\psi$ is trivial when restricted to 
$K_{p}\cap N_{p}$ for all $p$ not in $S$.)
Thus, in order to show that the Poincar\' e series is generic for some level $m$, 
it remains to show that the integral on the right is non-zero for some
matrix coefficient and some $g_{\infty}$ in $G(\mathbb R)$. This is 
done in the following section. 

\subsection{Some results of Wallach}
In this section $G=G(\mathbb R)$. Let $K$ be a maximal compact subgroup in 
$G$. Let $(\pi, H)$ be a discrete series representation on a Hilbert 
space $H$. Let $\langle v,w\rangle$ denote the inner product on $H$. 
Let $v$ be a non-zero vector in $H_{K}$, the space of 
$K$-finite vectors in $H$,  and consider the 
matrix coefficient $c_{v,w}(g)=\langle \pi(g)v, w\rangle$. It will be important for 
us that the function $c_{v,w}$ belongs to the (Harish-Chandra) 
Schwarz space $\mathcal C(G)$. 

Assume now that $\pi$ is a generic representation with respect to 
a regular unitary character $\psi$ of $N$.  The Whittaker functional $W_{\psi}$ 
 is not defined on $H$.  Instead, the Whittaker functional is defined 
on a space of smooth vectors $H^{\infty}$ and continuous with respect to a certain 
topology on $H^{\infty}$. 
Note that $H_{K}$, the space of 
$K$-finite vectors, is contained in $H^{\infty}$.  For every vector $v$ in $H_{K}$ 
we can define a generalized matrix coefficient 
$$ 
\ell_{\psi,v}(g)=W_{\psi} (\pi(g)v).
$$
Of course, $\ell_{\psi,v}(ng)=\psi(n)\ell_{\psi,v}(g)$ for every $N$. Moreover, 
the following important property of $\ell_{\psi,v}$ has been established by 
Wallach in \cite{Wal}, Section 15. 
The function $\ell_{\psi,v}$ belongs to the space of Schwarz functions 
 $\mathcal C(N\backslash G, \psi)$. This space is described using the 
 Iwasawa decomposition $G=NAK$. Here $A=\exp(\mathfrak a)$ where $\mathfrak a$
 is a maximal split Cartan subalgebra of the Lie algebra $\mathfrak g$ of $G$. 
 A  smooth function $f$ on $G$ belongs to 
 $\mathcal C(N\backslash G, \psi)$ if $f(ng)=\psi(n)f(g)$ and for every $X$ 
 in the enveloping algebra of $\mathfrak g$ and every positive integer $d$ there is a 
 constant $C$ such that for all $a$ in $A$ and $k$ in $K$
 $$ 
 |R_{X} f(ak)|\leq C \rho(a)(1+||\log(a)||)^{-d}
 $$
 where $R_{X}f$ is obtained by differentiating $f$ by $X$ from the right. Note that 
 this definition says, in essence, that the restriction of $f$ to $A$ is a usual 
 Schwarz function on $A$ multiplied by the modular character $\rho(a)$. The 
 Haar measure $dg$ on the group $G$ can be decomposed as 
 $$ 
 dg=dn ~ \rho^{-2}(a) da ~dk. 
 $$
 It follows that 
 the space $\mathcal C(N\backslash G, \psi)$ admits a natural 
 $G$-invariant inner product 
$$ 
(\varphi_{1},\varphi_{2})=\int_{AK} \varphi_{1}(ak)
\overline{\varphi_{2}(ak)} ~\rho^{-2}(a) da ~dk. 
$$
The absolute convergence of this integral is clear. 
In fact, as we shall need this observation in a moment, if $\varphi_{1}$ 
is in the Schwarz space,  
so $\varphi_{1}(ak)\leq C_{2,d}\rho(a)(1+||\log(a)||)^{-d}$ for any $d$, 
 and $\varphi_{2}(ak)\leq C_{2} \rho(a)$ then 
the integral is still absolutely convergent.  Indeed, up to a non-zero factor, the 
integral is bounded by 
$$ 
\int_{A} (1+||\log(a)||)^{-d}~da 
$$
which is absolutely convergent for a sufficiently large $d$. 

The map $v\mapsto \ell_{\psi,v}$ from $H_{K}$ to 
$\mathcal C(N\backslash G, \psi)$ is an intertwining map preserving inner 
products.  In particular the matrix coefficient $c_{v,w}$ can be written as 
$$ 
c_{v,w}(g)=(R(g)\ell_{\psi,v}, \ell_{\psi,w})
$$
where $R$ denotes the action of $G$  on $\mathcal C(N\backslash G, \psi)$ 
by right translations.

\begin{prop} Let $\psi$ be a regular (generic) unitary character of $N$. Let 
$(\pi, H_{K})$ be a $\psi$-generic discrete series. 
For every $v\neq 0$ in $H_{K}$ there are $g$ and $g_{1}$ in $G$ such that 
$$ 
\int_{N}c_{v,\pi(g_{1})v}(ng)\psi(n) ~dn\neq 0.
$$

\end{prop}
\begin{proof} 
The proof is based on the following lemma:
\begin{lemma} \label{convergence}
Let $\alpha$ be a function in $\mathcal C$ and $\varphi$ a 
function in $\mathcal C(N\backslash G, \psi)$.  Then there exists a constant $C$ 
such that 
$$ 
 \int_{G}|\alpha(g)|\cdot |\varphi(g_{1}g)| ~ dg
  \leq C \rho(a_{1})
 $$ 
for every  $g_{1}=n_{1}a_{1}k_{1}$ in $G$.
\end{lemma}

We shall postpone the proof of this lemma in order to finish the proof of  
proposition, first.  If we take $\alpha=c_{v,v}$ and $\varphi=\ell_{\psi,v}$, then 
the lemma assures that the integral 
$$ 
\int_{N\backslash G}\int_{G}c_{v,v}(g)\ell_{\psi,v}(g_{1}g)
\overline{\ell_{\psi,v}(g_{1})}~dg ~dg_{1}
$$ 
converges absolutely. Reversing the order of integration, we can rewrite this 
integral as 
$$ 
\int_{G}c_{v,v}(g)\overline{(R(g)\ell_{\psi,v}, \ell_{\psi,v})}
=||c_{v,v}||^{2}_{L^{2}(G)}\neq 0
$$
since, as we have remarked before, $(R(g)\ell_{\psi,v}, \ell_{\psi,v})=c_{v,v}(g)$. 
By Fubini's theorem, it follows that for some $g_{1}$ in $G$,
$$ 
0\neq \int_{G} c_{v, v}(g) \ell_{\psi,v}(g_{1}g) ~dg.
$$
The substitution $g:=g_{1}^{-1}g$ gives 
$$
0\neq\int_{G} c_{v, \pi(g_{1})v}(g) \ell_{\psi,v}(g) ~dg.
$$
Since this is an absolutely convergent integral over $G$, it can be written 
as a double integral over $N\backslash G \times N$. Then, by Fubini's theorem,
 there exists $g$ in $N\backslash G$ such that 
$$ 
0\neq \ell_{\psi,v}(g) \int_{N}c_{v,\pi(g_{1})v}(ng) \psi(n) ~dn.
$$
(Here we used that $\ell_{\psi,v}(ng)=\psi(n)\ell_{\psi,v}(g)$.)

It  remains to prove Lemma \ref{convergence}. Substituting $g:=g^{-1}g$ and writing $g=nak$ 
the integral (in the statement of Lemma \ref{convergence}) can be written as 
$$ 
\int_{NAK}|\alpha(k_{1}^{-1}a_{1}^{-1}n_{1}^{-1}nak)|\cdot 
|\varphi(nak)| ~dn~\rho^{-2}(a) da ~dk.
$$
Note that $|\varphi(nak)|=|\varphi(ak)|$ since $\psi$ is unitary.
We can use a substitution 
$n:=n_{1}n$ to rewrite the integral as 
$$ 
\int_{NAK}|\alpha(k_{1}^{-1}a_{1}^{-1}nak)|\cdot
|\varphi(ak)| ~dn~\rho^{-2}(a) da ~dk.
$$
Next, substituting $n:=a_{1}na_{1}^{-1}$ (this change of variable in $N$ contributes a 
factor $\rho^{2}(a_{1})$), the integral further becomes 
$$ 
\rho^{2}(a_{1})\int_{NAK}|\alpha(k_{1}na_{1}^{-1}ak)|\cdot 
|\varphi(ak)| ~dn~\rho^{-2}(a) da ~dk.
$$
Since $\alpha$ is in $\mathcal C(G)$, there exists a constant $c$ such that 
(see \cite{War}, Theorem 8.5.2.1)
$$ 
\int_{N}|\alpha(k_{1}^{-1}na_{1}^{-1}ak)| ~dn \leq c\rho(a_{1}^{-1}a)
$$
for all $(k_{1},k)$ in $K\times K$. It follows that the integral
(in the statement of Lemma \ref{convergence}) is bounded by 
$$ 
c\rho(a_{1}) \int_{AK}\rho(a)^{-1}|\varphi(ak)| ~da ~dk \leq C \rho(a_{1}), 
$$
for some constant $C$, exactly what we wanted. Lemma \ref{convergence} is proved. 
\end{proof}

\vskip 5pt 

Of course, our discussion is valid in the case of $p$-adic fields, provided that 
for every positive integer $d$, there exists a constant $C$ such that 
$$ 
|\ell_{\psi,v}(nak)| \leq C \rho(a) (1+ ||\log(a)||)^{-d}. 
$$
This may not be known in general, but if the discrete series is supercuspidal, then 
$\ell_{\psi,v}$ is compactly supported, so the proposition holds in this case, as 
well. Summarizing, we have shown the following: 

\begin{theorem} \label{global}
Let $G$ be a simple, quasi-split algebraic group defined over 
$\mathbb Q$. Fix two finite and disjoint sets of places: $D$ containing $\infty$ 
and perhaps nothing else,  and
$S$ a non-empty set of primes such that $G$ is unramified at all primes $p$ not in $D\cup S$. 
(This means that $G(\mathbb Q_{p})$ contains a hyperspecial maximal 
subgroup.)
Let $\psi$ be a regular (generic) character of $N(\mathbb A)$
trivial on $N(\mathbb Q)$. 
Assume that we are given a $\psi$-generic integrable discrete series representation 
of $G(\mathbb R)$, and a $\psi$-generic supercuspidal representation $\pi_{q}$
for every $q$ in $D$. Then there exists a global $\psi$-generic cuspidal representation 
$\pi$ such that $\pi_{\infty}$ is the given integrable discrete series, 
$\pi_{q}$ is the given supercuspidal representation for every $q$ in $D$ and 
$\pi_{p}$ is unramified for every $p$ outside $D\cup S$. 
\end{theorem}

\section{Proof of Theorem \ref{main}}\label{proofs}

We first  reduce the proof of Theorem \ref{main} to the construction of certain self-dual cuspidal automorphic representations on $\GL_n(\A_{\Q})$. Then we carry out the construction combining Theorem \ref{global} with the results in \cite{CKPSS}.

But to begin with, to apply Theorem \ref{global} to construct generic cuspidal representations with a given integral integrable discrete series at the infinite place on certain orthogonal groups, we need 
a description of generic, integrable discrete series representations of  the real group $\SO(m,m+1)$.

\subsection{Generic discrete series of  $\SO(m,m+1)$}\label{generic}

The Lie group $\SO(m+1,m)$ has two connected components. 
Let $G_{0}$ be the connected component containing the identity and $K_{0}$ a 
 maximal compact subgroup of $G_{0}$. Note that 
 $$ 
 K_{0}\cong \SO(m+1)\times \SO(m). 
 $$
 The necessary and sufficient condition for $G_{0}$ to have 
 discrete series representations is that the rank of $G_{0}$ is equal to the 
 rank of $K_{0}$. This clearly holds here. 
We shall now describe discrete series representations of $G_{0}$ and 
specify which of them are $\psi$-generic for a choice of a regular (generic) character 
$\psi$ of $N(\mathbb R)$, the unipotent radical of a Borel subgroup.
 (The difference between generic discrete series for
 $\SO(m+1,m)$ and $G_{0}$ is easy to explain. Any two generic characters of 
 $N(\mathbb R)$ are conjugate by an element in $\SO(m+1,m)$, whereas there are two conjugacy 
 classes of generic characters for $G_{0}$. Any generic discrete series 
 representation of $\SO(m+1,m)$, 
 when restricted to $G_{0}$, breaks up as a sum of two discrete series representation
 of $G_{0}$, each generic with respect to precisely one of the two classes of characters.)

 Let $\mathfrak g$ be the real Lie algebra of $G_{0}$ and $\mathfrak k$ the 
 real Lie algebra of $K_{0}$. 
 Let $\mathfrak h$ be a maximal Cartan subalgebra of $\mathfrak g$ 
 contained in $\mathfrak k$.  Let $\Phi$ and $\Phi_{K}$ be the sets of 
 roots for the action of $\mathfrak h$ on $\mathfrak g$ and $\mathfrak k$, 
 respectively. The roots in $\Phi_{K}$ are called {\em compact} roots. 
  The root system $\Phi$ is of type $B_{m}$. 
 We can identify $i\mathfrak h^{\ast}\cong \mathbb R^{m}$. 
 Let $(\cdot | \cdot)$ be the usual inner product on $\mathbb R^{m}$
such that the standard basis $e_{i}$, $1\leq i \leq m$ is orthonormal. 
Then  
 $$ 
 \Phi=\{ \pm e_{i}\pm e_{j}, \text{ with } i\neq j \text{ and } 
 \pm e_{i} \text{ for all } i\}.
 $$
 
 The Langlands parameter \cite[\S3]{Lan} defining an $L$-packet of discrete series 
 representation of $\SO(m+1,m)$ is a 
 homomorphism $\sigma_{\infty}: W_{\mathbb R}\rightarrow \Sp_{2m}(\mathbb C)$
 described as follows. Recall that $W_\R$ is the non-split extension of $\Z/2\Z$ by $\C^\times$ given by 
$W_\R=\C^\times \cup j\C^\times$ where $j^2=-1$ and $jzj^{-1}= \overline z$. 
The representation $\sigma_{\infty}$ 
is a direct sum of 2-dimensional symplectic
  representations $\rho_{i,\infty}$ ($ 1 \leq i \leq m$)  
which, when restricted to $\C^\times$, are of the form 
${({z \over \overline{z}})}^{{1-2\lambda_{i}}\over 2} 
\oplus {({z \over \overline{z}})}^{-{{1-2\lambda_{i}} \over 2}}$, 
for some non-zero integers $\lambda_{i}$ such that $\lambda_{i}\neq \pm \lambda_{j}$
if $i\neq j$, and $\rho_m(j)$ is the matrix
$$\left( \begin{array}{cc}
     0  & 1 \\
    -1 & 0
\end{array} \right).$$
The infinitesimal character of all representations in the $L$-packet of $\sigma_{\infty}$ is 
 $$ 
 \lambda=(\lambda_{1},\lambda_{2}, \ldots, \lambda_{m}) \in i\mathfrak h^{\ast}.
 $$
 In fact, the $L$-packet consists of all discrete series representations with this 
 infinitesimal character. 
More precisely, for every non-singular and integral $\lambda$ in $i\mathfrak h^{\ast}$,
 there exists a discrete series representation $\pi_{\lambda}$ of $G_{0}$
  with the infinitesimal character $\lambda$. 
 Furthermore, $\pi_{\lambda}\cong \pi_{\lambda'}$ if and only if $\lambda$ 
 and $\lambda'$ are conjugated by $W_{K}$, the Weyl group of $\Phi_{K}$. Of 
 course, $\pi_{\lambda}$ and $\pi_{\lambda'}$ have the same infinitesimal character
 if and only if $\lambda$ and $\lambda'$ are conjugated by $W$, the Weyl group of 
 $\Phi$. In particular, the number of representations in the $L$-packet (for $G_{0}$)
 is equal to the index of $W_{K}$ in $W$. 
 The representation $\pi_{\lambda}$ is $\psi$-generic for some choice of a 
 regular character $\psi$ of $N(\mathbb R)$ if and only if all walls of the 
 Weyl chamber containing  $\lambda$ are defined by non-compact roots
 (see \cite[\S6]{Vo} and \cite{Kos}). The 
 existence of one such chamber, in fact precisely two up to the action of $W_{K}$, 
 can be shown as follows. Instead of fixing an embedding $\Phi_{K}\subseteq \Phi$ 
 we shall fix a Weyl chamber $\mathcal C$ 
 containing $\lambda$, and then look for ways how to 
 put $\Phi_{K}$ into $\Phi$ so that it misses the roots defining the Walls of $\mathcal C$.
We pick the Weyl
chamber $\mathcal C$ containing $\lambda$, so that 
$\lambda=(\lambda_{1}, \lambda_{2},\ldots , \lambda_{m})$ where
 $\lambda _{i}$ are positive integers such that 
$\lambda_{1}> \cdots > \lambda_{m}$. In particular, the walls of the 
Weyl chamber are given by 
$$ 
e_{1}-e_{2}, e_{2}-e_{3}, \ldots e_{m-1}-e_{m} \text{ and } e_{m}. 
$$
Then $\pi_{\lambda}$ is $\psi$-generic for some choice of $\psi$ 
if and only if these roots are not compact. Thus, we 
need to show that we can embed $\Phi_{K}$ into $\Phi$ so that it does not 
contain any of these roots. To this end, break up the set of indices $\{1,2, \ldots, m\}$
into a disjoint union $E\cup O$ where $E=\{m, m-2,  \ldots \}$ 
and $O=\{m-1, m-3, \ldots \}$. Then we can pick $\Phi_{K}$ so that it 
contains $\pm e_{i}\pm e_{j}$ where {\em both} $i$ and $j$ are either in $E$ or 
in $O$, and $\pm e_{i}$ with $i$ in $O$. With this choice of $\Phi_{K}$, the 
discrete series $\pi_{\lambda}$ is generic.  The other Weyl chamber without 
``compact'' walls is $-\mathcal C$. These two are not $W_{K}$-conjugated since 
$-1$ is not contained in $W_{K}$. We put 
$$ 
\pi_{\infty}=\Ind^{\SO(m+1,m)}_{G_{0}}\pi_{\lambda}. 
$$
This is the unique generic discrete series representation of $\SO(m+1,m)$ with 
the infinitesimal character $\lambda$. In order to make this representation a 
local component of a global automorphic representation, we need that its matrix 
coefficients are integrable, as well. 
 Integrability conditions on matrix coefficients are given as follows (see \cite{Mi}):

 \begin{prop} Let $W$ be the Weyl group of $\Phi$. 
 Fix a positive $W$-invariant inner product 
 $(\cdot|\cdot)$ on $i\mathfrak h^{\ast}$.
  The discrete series representation $\pi_{\lambda}$ has 
 integrable matrix coefficients if 
 $$ 
|(\lambda|\alpha)| >k(\alpha)= \frac{1}{4}\sum_{\beta\in \Phi} |(\alpha|\beta)|.
 $$
 for every non-compact root $\alpha$. 
 \end{prop}
 
 In practical terms this simply means that $\lambda$ is at a certain distance from 
 all walls corresponding to non-compact roots. 
 We can determine whether the discrete series  $\pi_{\lambda}$ 
is integrable or not since one easily computes that 
$$ 
k(\alpha) =\begin{cases} m-\frac{1}{2} \text{ if $\alpha$ is short} \\
2m \text{ if $\alpha$ is long.}
\end{cases}
$$
In particular, if $\lambda_{m}\geq m$ and 
$\lambda_{i}-\lambda_{i+1}> 2m$ for all $i=1, \ldots m-1$ then 
$\pi_{\lambda}$ has integrable matrix coefficients. 

\subsection{Reduction of Theorem \ref{main} to existence of certain  cuspidal automorphic representations of $\GL_{n}(\A_\Q)$}\label{desired}

Let $\Pi$  be a cuspidal automorphic representation  of 
$\GL_n(\A_\Q)$ which is unramified or supercuspidal at each finite place $v$ of $\Q$. 
There is attached to $\Pi_v$ a representation $\sigma(\Pi_v)\colon W_{\Q_v} \rightarrow \GL_n(\aQl)$.
This arises from the local Langlands correspondence of \cite{HT} (for finite places, for infinite places these are the results of Harish-Chandra and Langlands, see \cite{Lan}, \cite{Bo}),
and depends on choosing an isomorphism $\C \simeq \aQl$.

Let $\Pi$  be a cuspidal automorphic representation  of 
$\GL_n(\A_\Q)$ with the following properties:

\begin{itemize}

\item (a) $\Pi$ is self-dual, i.e., $\Pi^{\vee} \simeq \Pi$

\item (b) $\Pi_\infty$ has  a regular symplectic parameter $\sigma_{\infty}$ 
described in the Section \ref{generic}. Recall that $\sigma_{\infty}$ is a 
is a direct sum of 2-dimensional  representations $\rho_{i,\infty}$ ($ 1 \leq i \leq m$)  
which, when restricted to $\C^\times$, are of the form 
${({z \over \overline{z}})}^{{1-2\lambda_{i}}\over 2} 
\oplus {({z \over \overline{z}})}^{-{{1-2\lambda_{i}} \over 2}}$. We require that
 $\lambda_{i}$ be positive integers such that $\lambda_{m}\geq m$ and 
$\lambda_{i}-\lambda_{i+1}> 2m$ for all $i=1, \ldots m-1$. This technical 
condition on the $\lambda_i$'s assures us that $\Pi_{\infty}$ is a local 
lift of an integrable discrete series representation $\pi_{\infty}$ of $\SO(m+1,m)$. 
\item (c) $\Pi$ is unramified outside $\{\ell,q\}$, and $\sigma(\Pi_q)$ is isomorphic to the $\rho_q$ fixed above.
\end{itemize}

The results of  \cite{Kot}, \cite{Cl}, \cite{HT},
see Theorem 3.6 of \cite{Tay} (applied to a twist of $\Pi$ by the $1-n \over 2$
power of the norm character), ensure that there is a  continuous semisimple representation $\rho_{\Pi}:G_\Q \rightarrow \GL_n(\aQl)$ attached
to $\Pi$ such that  for the finite places $v \neq \ell$, $\rho_{\Pi}|_{D_v} \simeq \sigma(\Pi_v) \otimes | \ \ |^{{1-n} \over 2}$. 
Here $| \ \ |^{1 \over 2}: G_{\Q_q} \rightarrow \aQl$ is  the unramified character of $\Q_q^\times$ that takes $q \rightarrow \sqrt{q}$ ($\sqrt{q}$ positive). For any integer $r$, we may also analogously  define a character $| \ \ |^r$ of $G_\Q$ with values in $\aQl^*$ which is the $r$th power of the $\ell$-adic cyclotomic character.

From the fact that $\Pi$ is self-dual
we see by \v Cebotarev density
that $\rho_{\Pi}^\vee \simeq \rho_{\Pi} | \ \ |^{n-1}$ and thus $\rho_{\Pi}$ acts either by orthogonal or symplectic similitudes on $\aQl^n$ with similitude factor $ | \ \ |^{n-1}$. 
Although it is possible for an irreducible representation to act by both orthogonal and
symplectic similitudes, this is not possible if the factors of similitude are the same.
As $\rho_{\Pi}|_{D_q} \simeq \rho_q \otimes | \ \ |^{{1-n} \over 2}$, and $\rho_q$ is an irreducible symplectic representation, it follows  
that $\rho_\Pi$ is irreducible, and that the self-duality of
$\rho_\Pi$ with similitude factor $|\ \ |^{n-1}$ is symplectic.
Therefore, the image of $\rho_{\Pi}$ may be conjugated to land inside $\GSp_n(\aQl)$, and in fact by the compactness of $G_\Q$, inside $\GSp_n(\aZl)$.

We consider the reduction mod $\ell$ of $\rho_{\Pi}$, and denote the resulting representation by $\rhobar:G_\Q \rightarrow \GSp_n(\bar\F_\ell)$, and note that its determinant is valued in $\F_\ell^\times$. Let $\Gamma$ denote ${\rm im}(\rhobar)$. Then we see that $\Gamma$ satisfies the conditions of Theorem \ref{crucial}, 
by construction, namely the choice of $q$
 and the parameter $\rho_q$. We expand on this. We see that any 
subgroup of $\Gamma$ of index $\leq N_d$ cuts out an extension $K$ of $\Q$ of degree 
$\leq N_d$ that is unramified outside $\{\ell,q,\infty\}$.  In fact as the image of $\rho_q(I_q)$ is of order $p$ and $p>N_d$ we see that  $K$ is unramified at $q$, and thus $K$ is unramified outside 
$\{\ell,\infty\}$. Thus by choice of $q$, it splits in $K$. Furthermore,  by construction ${\rm im}(\rhobar(D_q))$ is a group of type $(n,p)$ (note that by choice $p \neq \ell$) of Section \ref{groups}, and it is contained in $\Gamma^{N_d}$.

Thus Theorem \ref{crucial} implies that after conjugation by an element in $\GL_n(\bar\F_\ell)$ we may conclude 
that $\Gamma$ contains $\Sp_n(\F_{\ell^k})$  for some integer $k$ and is contained in its normalizer.  Thus 
by Corollary \ref{tobeused} we know that the image of $\Gamma$ in $\PGL_n(\bar \F_\ell)$ is isomorphic  to $\PSp_n(\F_{\ell^k})$ or $\GSp_n(\F_{\ell^k})/\F_{\ell^k}^\times$. 
As the order of this group is divisible by $p$ (as the order of $\rhobar(I_q)$ is $p$), it follows using Lemma \ref{pq}(5) that $k$ is divisible by $t$.

\subsection{Construction of certain cuspidal automorphic  representations of $\GL_n(\A_\Q)$}\label{last}
In order to construct the $\Pi$ of the previous section we construct generic cuspidal automorphic representations $\pi$ of the split  $\SO_{2m+1}(\A_\Q)$ using Theorem \ref{global}  and
lift them to $\GL_{2m}(\A_\Q)$ using the results of \cite{CKPSS}, \cite{JS1}, \cite{JS2}. We use the terminology of these papers below. 
 There is forthcoming work of Chenevier and Clozel which uses related, but more elaborate, constructions to improve the results in \cite{Che}, which Chenevier had mentioned to the first named author. The relevance of  $\SO_{2m+1}$ to our work is that the connected component of its $\rm L$-group is $\Sp_{2m}$.

We consider the split group $\SO_{2m+1}$ of rank $m$ defined over $\Q$ (defined by the form $\Sigma_{i=1}^nx_ix_{n+i} + x_{2n+1}^2$), and consider $\SO_{2m+1}(\Q_v)$ for each place $v$ of $\Q$, and
$\SO_{2m+1}(\A_\Q)$. We note that the notion of genericity  for these groups is independent of choice of (local or global) Whittaker character $\psi$, and thus we call the $\psi$-generic forms, or $\psi$-generic local representations,  of Section \ref{Poincare}  simply generic.
 
 We need the following theorem which is a combination of the work of 
 \cite{CKPSS} and \cite{JS2}: see \cite[Theorem 7.1]{CKPSS} and 
 \cite[Theorem E]{JS2}.

\begin{theorem}\label{lifting}
  There is a lifting from equivalence classes of irreducible generic cuspidal automorphic  
  representations  of  $\SO_{2m+1}(\A_\Q)$ to equivalence classes of irreducible 
  automorphic representations of $\GL_{2m}(\A_\Q)$ such that this lifting is functorial at 
  all places. Further a cuspidal automorphic representation $\Pi$ of  $\GL_{2m}(\A_\Q)$ which is in the image of this lift is self-dual (and $L(s,\Lambda^2,\Pi)$ has a simple pole at $s=1$).
\end{theorem}

We refer to the cited papers for the exact notion of functoriality used, 
but will spell it out in the cases used below. 

In order to construct the generic cuspidal representation $\pi$ we need to specify 
what we want at the local places. 
We start with the following theorem of Jiang-Soudry: \cite[Theorem  6.4]{JS1} 
and \cite[Theorem 2.1]{JS2}.

\begin{theorem}\label{js}
Let $q$ be finite prime of $\Q$. There is a bijection between irreducible generic discrete series representations of $\SO_{2m+1}(\Q_q)$ and irreducible generic representations of $\GL_{2m}(\Q_q)$ with Langlands parameter of the form $\sigma=\Sigma \sigma_i$ with $\sigma_i$ irreducible symplectic representations of $WD_{\Q_q}$ which are pairwise non-isomorphic.
\end{theorem}

 Thus in particular there is a generic supercuspidal representation $\pi_q$ of 
 $\SO_{2m+1}(\Q_q)$ that corresponds to 
the Langlands parameter $\rho_q$ (and thus to a supercuspidal representation of $\GL_{2m}(\Q_q)$
with this parameter).
This correspondence is also known at the Archimedean places as recalled in 
Section \ref{generic}.  From this we deduce there is a generic, integrable discrete series representation $\pi_\infty$ on $\SO_{2m+1}(\R)$ which corresponds 
(under the correspondence of \cite[Section 5.1]{CKPSS}) to the representation $\Pi_\infty$ fixed in Section \ref{desired} with Langlands parameter $\sigma_{\infty}$.

By Theorem \ref{global} (with $D=\{\infty,q\}$ and $S=\{\ell\}$) 
there exists a generic cuspidal automorphic representation 
$\pi$ on $\SO_{2m+1}(\A_\Q)$ such that:

\begin{itemize}
\item Under the Jiang-Soudry correspondence  
of Theorem \ref{js}, $\pi_q$ has parameter $\rho_q$.

\item $\pi$ is unramified outside $\{\ell,q\}$

\item  $\pi_{\infty}$ is a  generic integrable discrete series with Langlands parameter
$\sigma_{\infty}$.

\end{itemize}

Using Theorem \ref{lifting} we can transfer $\pi$ to $\Pi$ to get an irreducible 
automorphic representation $\Pi$ on $\GL_{2m}(\A_\Q)$ such that
\begin{itemize}

\item  $\Pi_\infty$ has the regular algebraic parameter $\sigma_{\infty}$,  $\Pi$ is unramified outside $\{\ell,q\}$, and $\sigma(\Pi_q)\simeq 
\rho_q$. (this for us is the implication of the {\it functorial at all places} assertion 
in Theorem \ref{lifting}).

\item $\Pi$ is cuspidal (as $\Pi_q$ is supercuspidal) and self-dual.

\end{itemize}

\noindent{\bf Remarks:}

\begin{itemize}

\item  It seems difficult to directly
construct self-dual representations of $\GL_{n}(\A_\Q)$ interpolating finitely many specified self-dual supercuspidal representations at finitely many places. As pointed out in \cite{PR1}, one of the difficulties is that an obstruction to this is that the corresponding local Langlands parameters  should either be all symplectic or all orthogonal, and proofs using the trace formula might not see 
this obstruction.  This is why we first construct $\pi$ on $\SO_{2m+1}(\A_\Q)$ and then transfer it to  $\GL_{2m}(\A_\Q)$ using the results of \cite{CKPSS}.

\item The case $n=2$ corresponds to the result of \cite{Wiese}. In that case
the lifting proved in \cite{CKPSS} is trivial: it is the lifting of cuspidal automorphic representations of $\PGL_2(\A_\Q)$ to cuspidal automorphic representations of $\GL_2(\A_\Q)$ with trivial central character.

\item Curiously enough as we lack control of the field of definition of the $\rhobar$ we get, we don't see using this method  how to realize 
$\PGL_2(\F_{\ell^k})$ as  a  Galois group over $\Q$ for infinitely many $k$.  The limitations of our method do not allow us to prove that given an integer $t>1$ there are infinitely many $k$ prime to $t$ (or even one such $k$) 
such that $\PSp_n(\F_{\ell^k})$
appears as a Galois group over $\Q$.

\item  In the $n=2$ case, we may prove the result of \cite{Wiese} for $\ell>2$
by imposing in addition to large dihedral ramification at a prime $q$, also $A_4/S_4$-type ramification at another prime (necessarily 2!). This works for $\ell>2$ to ensure that we get some large image representations, but does not work for $\ell=2$.  This is because by our methods it is not possible to ensure that a non-trivial  unipotent is in the image of the mod $\ell$ Galois representation being considered. A similar remark applies for higher dimensions $n$. 
Further it seems of interest to us to force large images  of global Galois representations by dint of properties of the representation at a {\it single}  prime $q$.

\item In an earlier version of the paper
 \\ (see {\tt http://front.math.ucdavis.edu/math.NT/0610860})  
it had been erroneously asserted
that the existence of generic cuspidal forms as in Theorem \ref{global} follows from the literature, in particular the methods of \cite{PS}. But it turns out that the methods of \cite{PS} using the relative trace formula are not able to 
prove results like Theorem \ref{global} where one of the local representations
sought to be interpolated into a global generic representation is a generic discrete series representation of a real group.

\end{itemize}

\section{Zariski density}\label{extra}

We conclude with a group-theoretic 
proposition which shows that if $t\gg 0$, the $\ell$-adic representations
$\rho_{\Pi}:G_\Q \rightarrow \GL_n(\aQl)$ constructed in Section \ref{proofs} have Zariski-dense image in 
$\GSp_n$.  

Before stating it, we first prove a lemma that we will need to prove it:

\begin{lemma}
\label{subquot}
Given a positive even integer $n$ and a prime $\ell$ 
there exists a constant $M$ such that for all $m>M$,
and all almost simple algebraic groups $G/\bar\F_\ell$ of rank $< n/2$,
the finite simple group $\PSp_n(\F_{\ell^m})$ is not a subquotient of $G(\bar\F_\ell)$.
\end{lemma}

\begin{proof}
Up to isomorphism there are only finitely many possibilities for $G$, so we may pick one.
Let  $r<n/2$ denote the rank of $G$, and $e_1<e_2<\cdots<e_r$ the exponents.
Let $p>e_r$ be any prime and $\F$ a finite field in characteristic $\ell$ 
such that $G$ is defined and split over $\F$ and
$p$ divides the order of $\F^\times$.  Let $T$ be an $\F$-split maximal torus of $G$.
We have
$$\ord_p |T(\F)| = r\ord_p (|\F|-1) = \sum_{i=1}^r \ord_p (|\F|^{e_i}-1) = \ord_p |G(\F)|,$$
so any $p$-Sylow subgroup of $T(\F)$ is a $p$-Sylow subgroup of $G(\F)$.  It follows that
every $p$-Sylow of $G(\F)$ is abelian and generated by $\le r$ elements, and these
properties are inherited by any finite $p$-subgroup of $G(\F)$ and therefore (letting $\F$ grow)
of $G(\bar\F_\ell)$.  It follows that no finite subgroup of $G(\bar\F_\ell)$ has a subquotient isomorphic to $(\Z/p\Z)^{n/2}$.  By
Lemma~\ref{cyclotomic-primes}, for $m$ sufficiently large, $\ell^m-1$ has a prime divisor 
 $p>e_r$, so
$\PSp_n(\F_{\ell^m})$ has a subgroup isomorphic to $(\Z/p\Z)^{n/2}$.
It  cannot, therefore, be a subquotient of $G(\bar\F_\ell)$. 

\end{proof}

Let $\Gamma = \rho_\Pi(G_{\Q})$.  The image of $\Gamma$
lies in $\GL_n(K)$ for some $\ell$-adic field $K$.  Since $\rho_\Pi$ has positive weight, in order to
prove that the Zariski-closure of $\Gamma$ is $\GSp_n$, it suffices to prove that the closure
contains $\Sp_n$.  This follows from the following proposition:

\begin{prop}
Let $K$ be a finite extension of $\Q_\ell$ with residue field $k$
and $\Gamma$ denote a compact subgroup of $\GSp_n(K)\subset \GL_n(K)$.
Suppose that  some quotient of $\Gamma$ is isomorphic to $\PSp_n(\F_{\ell^m})$.  
If $m$ is sufficiently large, then the Zariski-closure of $\Gamma$ contains
$\Sp_n$.
\end{prop}

\begin{proof}
Let $G$ denote the Zariski-closure of $\Gamma$ in $\GL_n$.  
Let $G^\circ$ denote the identity component of $G$ and
$H$ the quotient $G/G^\circ$.   There is a version of Jordan's theorem for algebraic groups
which asserts that $H$ contains a normal abelian subgroup of index $\le J(n)$.   
This seems to be well-known but as we cannot locate a reference, 
we sketch a proof here.

It suffices to prove that for some finite abelian extension $\tilde H$ of $H$, the homomorphism
$\tilde H\to H$ lifts to $\tilde H\to G(\C)$.  Indeed Jordan's theorem for finite subgroups of
$\GL_n(\C)$ then applies to $\tilde H$, and the image of any normal abelian subgroup of $\tilde H$
is a normal abelian subgroup of $H$.   Lifting by stages, it suffices to prove this first
in the case the $G^\circ$ is adjoint semisimple, next in the case that $G^\circ$ is diagonal, and last
in the case that $G$ is commutative and unipotent.  
For the first case, we note that the center of $G(\C)$ is trivial, so every extension of 
$H$ by $G(\C)$ is a semidirect product.  For the second, we note that 
$H^2(H,D(\C))$ is annihilated by $|H|$, and therefore lies in the image of
$H^2(H,D(\C)[|H|])$.  Thus, every cohomology class in $H^2(H,D(\C))$
can be trivialized by pullback to a finite abelian extension $\tilde H$ of $H$.
For the third, we note that $H^2(H,V) = 0$
for every complex representation $V$ of $H$, so there is no obstruction to lifting.

Now, if $0\to G_1\to G_2\to G_3\to 0$ is any short exact sequence of groups and $G_2$
admits a surjective homomorphism to a finite simple group $\Delta$, then $G_1$ maps to a normal
subgroup of $\Delta$; thus either $G_1$ or $G_3$ maps onto $\Delta$.   
Setting $\Delta =  \PSp_n(\F_{\ell^k})$ and assuming $|\Delta| > J(n)$, we see that
the component group $H$ cannot map to $\Delta$, and therefore
$G^\circ(K)\cap \Gamma$ must.  Without loss of generality, therefore, we may assume that
$G$ is connected.  
If $R$ denotes the radical of $G$, then $R(K)\cap\Gamma$ is a normal 
prosolvable subgroup of $\Gamma$, so its image in $\Delta$ is trivial.  It follows that
there exists a semisimple quotient $G_s$ of $G$ such that $G_s(K)$ contains a compact subgroup 
$\Gamma_s$ which admits a surjective homomorphism to $\Delta$.  
Replacing $K$ with a finite extension $L$, we may assume that
$\Gamma_s$ stabilizes a hyperspecial vertex of the building of $G_s$ over $L$ 
(\cite[Prop.~8]{Se2}, \cite[Lemma~2.4]{Lar}).
It follows that there exists a smooth group scheme $\cG_s$ over the ring of integers $\cO_L$ of $L$
with connected semisimple fibers such that $\Gamma_s\subset \cG_s(\cO_L)$
and the generic fiber of $\cG_s$ is isomorphic to $G_s$.  The kernel of the reduction map on $\Gamma_s$ is a normal pro-$\ell$-group of $\Gamma_s$ whose image in
$\Delta$ must again be trivial.  We conclude that the image of $\Gamma_s$ under
the reduction map admits a surjective homomorphism to $\Delta$.
Let $G_s^\ell$ denote the special fiber of $\cG_s$.  It is connected and semisimple, with the same Dynkin diagram as $G_s$.  Moreover $G_s^\ell(\bar\F_\ell)$ contains a subgroup 
which maps onto $\Delta$.

We assume that $G_s$, or equivalently $G_s^\ell$, is not symplectic of rank $n/2$.
If the rank of $G_s$ is $n/2$ but $G_s\neq \Sp_n$, 
then by the classification of equal rank subgroups of $\Sp_n$,
$G_s$ fails to be almost simple.  In this case, we can replace $G_s^\ell$ by an almost simple subquotient,
whose rank is strictly less than $n/2$.  In any case, as long as $G_s\neq \Sp_n$,
we can find $G_s^\ell$ with rank less than
$n/2$ such that $G_s^\ell(\bar\F_\ell)$ contains a finite subgroup which maps onto 
$\Delta = \PSp_n(\F_{\ell^m})$.  By Lemma~\ref{subquot}, this cannot happen for $m\gg 0$.

\end{proof}

\section{Acknowledgements}

  We are much indebted to  Dragan Mili\v ci\'c, Dipendra Prasad, Peter Trapa 
  and Nolan Wallach for their assistance, especially with real groups.
   We thank Michael Dettweiler for some helpful correspondence, and Don Blaius for helpful feedback on the manuscript.

\end{document}